\documentclass[a4paper,12pt]{amsart}
\usepackage{amssymb} \usepackage{amsfonts}
\usepackage{amsmath,amsthm}
\usepackage{tikz,relsize}
\usepackage[cmtip,arrow]{xy}
\usepackage{pb-diagram,pb-xy}

\newcommand{\Z}{{\mathbb Z}}

\newcommand{\R}{{\mathbb R}}

\newtheorem{lm}{Lemma}[section]

\newtheorem{theo}{Theorem}[section]

\theoremstyle{definition}

\newtheorem{ex}{Example}[section]
\begin{document}
	
\title[Note about Stiefel-Whitney classes on real Bott manifolds ]{Note about Stiefel-Whitney classes on real Bott manifolds}

\author{A. G\c{a}sior}

\begin{abstract}  Real Bott manifolds is a class of flat manifolds with holonomy group $\mathbb Z_2^k$
of diagonal type. In this paper we want to show how we can compute even Stiefel - Whitney classes on real Bott manifolds. This paper is an answer to the question of professor Masuda if is it possible to extend \cite{GS2} and compute any Stiefel-Whitney classes for real Bott manifolds.  It also extends results of \cite{GS1}.
\end{abstract}
\keywords{Real Bott manifolds, Stiefel-Whitney class
\newline Author is supported by the Polish National Science Center grant DEC-2017/01/X/ST1/00062}

\maketitle

\hskip5mm
\section{Introduction}
Let $M_n$ be a flat manifold of dimension $n$, i.e. a compact connected Riemannian manifold without boundary
with zero sectional curvature. From the theorem of Bieberbach (\cite{Ch}, \cite{S3})
the fundamental group
$\pi_{1}(M_{n}) = \Gamma$ determines a short exact sequence:
\begin{equation}\label{ses}
0 \rightarrow \Z^{n} \rightarrow \Gamma \stackrel{p}\rightarrow
G \rightarrow 0,
\end{equation}
where
$\Z^{n}$ is a maximal torsion free abelian subgroup of rank $n$ and
$G$ is a finite group which
is isomorphic to the holonomy group of $M_{n}.$
The universal covering of $M_{n}$ is the Euclidean space $\R^{n}$
and hence $\Gamma$
is isomorphic to a discrete cocompact subgroup
of the isometry group $\operatorname{Isom}(\R^{n}) = \operatorname{O}(n)\times\R^n = E(n).$ In that case $p:\Gamma\to G$ is a projection on the first component of the semidirect product $O(n)\ltimes \mathbb R^n$ and $\pi_1(M_n)=\Gamma$ is a subgroup of $O(n)\ltimes \mathbb R^n$.
Conversely, given a short exact sequence of the form (\ref{ses}), it is known that
the group $\Gamma$ is (isomorphic to) the fundamental group of a flat manifold.
In this case $\Gamma$ is called a Bieberbach group.
We can define a holonomy representation $\phi:G\to \operatorname{GL}(n,\Z)$ by the formula:
\begin{equation}\label{holonomyrep}
\phi(g)(e) = \tilde{g}e(\tilde{g})^{-1},
\end{equation}
for all $e\in\mathbb Z^n, g\in G$ and where $p(\tilde{g})=g.$ In this article we shall consider Bieberbach groups of rank $n$ with holonomy group  $\Z_{2}^{k}$, $1\leq k\leq n-1$,
and  $\phi(\Z_{2}^{k})\subset D\subset \operatorname{GL}(n,\Z)$.
Here $D$ is the group of matrices with $\pm1$ on the diagonal.

The main result is the formula for even Stiefel-Whitney classes for real Bott manifolds. This formula is generalization of the one from our previous paper (\cite{GS2}, Lemma 2.1). It was suggested to us by M. Masuda.
The author thanks Andrzej Szczepa\'{n}ski for discussion.

\section{Stiefel-Whitney classes for real Bott manifolds}
Let
\begin{equation}\label{tower}
M_{n}\stackrel{\R P^1}\to M_{n-1}\stackrel{\R P^1}\to...\stackrel{\R P^1}\to M_{1}\stackrel{\R P^1}\to M_0 =
\{ \bullet\}
\end{equation}
be a sequence of real projective bundles such that $M_i\to M_{i-1}$, $i=1,2,\ldots,n$, is a projective
bundle of a Whitney sum of a real line bundle $L_{i-1}$ and the trivial line bundle over $M_{i-1}$.
The sequence (\ref{tower}) is called the real Bott tower and the top manifold $M_n$ is called the real
Bott manifold, \cite{CMO}.

Let $\gamma_i$ be the canonical line bundle over $M_i$ and we set $x_i = w_1(\gamma_i)$ ($w_1$ is the first Stiefel-Whitney class).
Since $H^1(M_{i-1},\Z_2)$ is additively generated by $x_1,x_2,..,x_{i-1}$ and $L_{i-1}$
is a line bundle over $M_{i-1},$ we can uniquely write
\begin{equation}\label{w1}
w_1(L_{i-1}) = \sum_{l=1}^{i-1} a_{li}x_l
\end{equation}
where $a_{li}\in \Z_2$ and $i = 2,3,...,n.$

From the above we obtain the matrix $A = [a_{li}]$ which is an $n\times n$ strictly upper triangular matrix
whose diagonal entries are $0$ and remaining entries are either
$0$ or $1.$
 One can observe (see \cite{KM}) that the tower (\ref{tower}) is completly determined by the matrix $A$ and therefore we may denote the real Bott manifold $M_n$ by $M_n(A)$.
From \cite[Lemma 3.1]{KM} we can consider $M_n(A)$ as the orbit space  $M_n(A) = \R^n/\Gamma(A),$
where $\Gamma(A)\subset E(n)$ is generated by elements
$$	s_{i} = \left(\operatorname{diag}\left[	1,\ldots,(-1)^{a_{i,i+1}},\ldots,(-1)^{a_{i,n}}\right],\left(0,\ldots,0,\frac{1}{2},0,\ldots,0\right)^T\right),$$
 where $(-1)^{a_{i,i+1}}$ is in the $(i+1, i+1)$ position and $\frac{1}{2}$ is the $i$th coordinate of the last column, $i = 1,2,...,n-1.$
$s_{n} = \left(I,\left(0,0,...,0,\frac12\right)\right)\in E(n).$
From \cite[Lemma 3.2, 3.3]{KM} $s_{1}^{2},s_{2}^{2},...,s_{n}^{2}$ commute with each
other and generate a free abelian subgroup $\Z^n.$ In other words $M_n(A)$ is a flat manifold with holonomy group $\mathbb Z_2^k$ of
diagonal type. Here $k$ is a number of non zero rows of a matrix $A$.

We have the following two lemmas.

\begin{lm}[\cite{KM}, Lemma 2.1]\label{lemma11}
The cohomology ring $H^*(M_n(A),\mathbb Z_2)$ is generated  by degree one elements $x_1,\ldots,x_n$ as a graded ring with $n$ relations
$$x_j^2=x_j\sum_{i=1}^na_{ij}x_i,$$
for $j=1,\ldots,n$.
\end{lm}
\begin{lm}[\cite{KM}, Lemma 2.2]\label{lemma12}
The real Bott manifold $M_n(A)$ is orientable if and only if the sum of entries is $0 (\operatorname{mod}2)$ for each row of the matrix $A$.
\end{lm}

The $k$th Stiefel-Whitney class \cite[ page 3, (2.1) ]{LS}  is given by the formula
\begin{equation}
w_k(M(A)) = (B(p))^{\ast}\sigma_{k}(y_1,y_2,...,y_{n})\in H^{k}(M(A);\Z_2) ,
\end{equation}
where $\sigma_k$ is the $k-$th elementary symmetric function, $B(p)$ is a map
induced by $p$ on the classification space and
\begin{equation}
y_i : = w_1(L_{i-1})\label{y}\end{equation}
 for $i=2,3,\ldots,n$.

Follow \cite{CMM}, if we consider $H^*(M_j(A),\mathbb Z)$ as a subring of $H^*(M_n(A),\mathbb Z)$ through the projection in (\ref{tower}), we see that
\begin{equation}\begin{aligned}
&H^*(M_n(A),\mathbb Z)\\
&=\mathbb Z[x_1,\ldots,x_{n}]/ \left(x_j^2-x_j\sum_{i=1}^na_{ij}x_i:j=1,2,\ldots,n\right).
\end{aligned}\end{equation}
From the above we get
\begin{lm}\cite{CMM}
Let $k$ be positive integer less or equal to $\frac n2$. The the set
$$\{x_{i_1}x_{i_2}\ldots x_{i_{2k}}:1\leq i_1<i_2<\ldots< i_{{2k}}\leq n\}$$
is an additive basis of $H^{2k}(M_n(A),\mathbb Z_2).$
\end{lm}

Let $A_{i_1i_2\ldots i_{2k}}$ denotes the $(n\times n)$ matrix consisting of $i_1,i_2,\ldots, i_{2k}$ rows of matrix $A$. Then non zero entries are only in $i_1,i_2,\ldots, i_{2k}$ rows of the matrix $A_{i_1i_2\ldots i_{2k}}$ and we have the following main result.

\begin{theo}\label{lemat2}
    Let $A$ be an $(n\times n)$ the Bott matrix.
    Then,
    $$w_{2k}(M_n(A))=\sum_{1\leq i_1<i_2<\ldots<i_{2k}\leq n}w_{2k}(M_n(A_{i_1i_2\ldots i_{2k}})).$$
\end{theo}

\noindent
{\bf Proof.}

From (\cite{CMM} Lemma 2.1) we have that the $2k$ cohomology group of $H^{2k}(M_n(A),\mathbb Z_2)$ has a basis
$${\mathcal{B}}=\{x_{i_1}x_{i_2}\ldots x_{i_{2k}}:1\leq i_1<i_2<\ldots<i_{2k}\leq n\}.$$
Moreover, also from Lemma \ref{lemma11} $x_j^2$ can be expressed by a linear combination of $x_kx_j$ for $k<j$. Note that this combination always contains an $x_j-$term. Hence, we get that $w_{2k}(M_n(A))$ is a sum of linear elements
$$w_{2k}(M_n(A))=\sum_{1\leq i_1<i_2<\ldots<i_{2k}\leq n}x_{i_1}x_{i_2}\ldots x_{i_{2k}}.$$
Each term $x_{i_1}x_{i_2}\ldots x_{i_{2k}}$ of this sum is an element from basis ${\mathcal{B}}$ and it is equal to the $2k$ Stiefel-Whitney class of the real Bott manifold $M_n(A_{i_1i_2\ldots i_{2k}})$, so we get
$$w_{2k}(M_n(A))=\sum_{1\leq i_1<i_2<\ldots<i_{2k}\leq n}w_{2k}(M_n(A_{i_1i_2\ldots i_{2k}})).$$
Thus, the $2k$th Stiefel-Whitney class of the real Bott manifold $M_n(A)$
 is equal to the sum of $2k$th Stiefel-Whitney classes of elementary components $M_n(A_{i_1i_2\ldots i_{2k}})$, $1\leq i_1<i_2<\ldots<i_{2k}\leq n$.

\hskip 142mm $\Box$

At the end of the paper we give an example.
\begin{ex}
For
$$A=\left[\begin{matrix}0&1&1&0&0&0&0\\0&0&1&1&0&0&0\\0&0&0&1&1&0&0\\0&0&0&0&1&1&0\\0&0&0&0&0&1&1\\0&0&0&0&0&0&0\\0&0&0&0&0&0&0\end{matrix}\right]$$
we get $w_4(M(A))=x_2x_3x_4x_5+x_1x_3x_4x_5+x_1x_2x_3x_5+x_1x_2x_3x_4.$
For the matrix $A$ we have the following
$$\begin{aligned}A_{1234}&=\left[\begin{matrix}0&1&1&0&0&0&0\\0&0&1&1&0&0&0\\0&0&0&1&1&0&0\\0&0&0&0&1&1&0\\0&0&0&0&0&0&0\\0&0&0&0&0&0&0\\0&0&0&0&0&0&0\end{matrix}\right],\; w_4\left(M\left(A_{1234}\right)\right)=x_1x_2x_3x_4,\\
A_{1235}&=\left[\begin{matrix}0&1&1&0&0&0&0\\0&0&1&1&0&0&0\\0&0&0&1&1&0&0\\0&0&0&0&0&0&0\\0&0&0&0&0&1&1\\0&0&0&0&0&0&0\\0&0&0&0&0&0&0\end{matrix}\right], \; w_4\left(M\left(A_{1235}\right)\right)=x_1x_2x_3x_5,\\
A_{1245}&=\left[\begin{matrix}0&1&1&0&0&0&0\\0&0&1&1&0&0&0\\0&0&0&0&0&0&0\\0&0&0&0&1&1&0\\0&0&0&0&0&1&1\\0&0&0&0&0&0&0\\0&0&0&0&0&0&0\end{matrix}\right], \; w_4\left(M\left(A_{1245}\right)\right)=0,\end{aligned}$$
$$\begin{aligned}A_{1345}&=\left[\begin{matrix}0&1&1&0&0&0&0\\0&0&0&0&0&0&0\\0&0&0&1&1&0&0\\0&0&0&0&1&1&0\\0&0&0&0&0&1&1\\0&0&0&0&0&0&0\\0&0&0&0&0&0&0\end{matrix}\right], \; w_4\left(M\left(A_{1345}\right)\right)=x_1x_3x_4x_5,\\
A_{2345}&=\left[\begin{matrix}0&0&0&0&0&0&0\\0&0&1&1&0&0&0\\0&0&0&1&1&0&0\\0&0&0&0&1&1&0\\0&0&0&0&0&1&1\\0&0&0&0&0&0&0\\0&0&0&0&0&0&0\end{matrix}\right], \; w_4\left(M\left(A_{2345}\right)\right)=x_2x_3x_4x_5.
\end{aligned}$$
So we have
$$\begin{aligned}
&w_4(M(A_{1234}))+w_4(M(A_{1235}))+w_4\left(M\left(A_{1245}\right)\right)+w_4(M(A_{1345}))+w_4(M(A_{2345}))\\&=x_1x_2x_3x_4+x_1x_2x_3x_5+x_1x_3x_4x_5+x_2x_3x_4x_5=w_4(M(A)).
\end{aligned}$$
\end{ex}

\vskip 2mm
\noindent
 Maria Curie-Sk{\l}odowska University,\\
  Institute of Mathematics\\
pl. Marii Curie-Sk{\l}odowskiej 1\\
20-031 Lublin, Poland\\
E-mail: anna.gasior@poczta.umcs.lublin.pl

\end{document}